\documentclass[a4paper,11pt]{article}
\usepackage[hidelinks]{hyperref}
\hypersetup{
    colorlinks=true,
    linktoc=all,
    linkcolor=red,     
   citecolor=blue,     
}
\usepackage{tikz}
\usepackage{subfigure}
\usepackage{graphicx}
\usepackage{amssymb}
\usepackage{latexsym, bm}
\usepackage{multicol}
\usepackage{indentfirst}
\usepackage{amssymb,amsfonts}
\usepackage{amsmath}
\usepackage{hyperref}

\usepackage{diagbox}
\usepackage[utf8]{inputenc}
\usepackage{authblk}
\usepackage{setspace}
\newcommand{\dd}{\delta}
\newcommand{\de}[1]{\delta_{#1}}
\def\twr{\mbox{\rm tw}}

\newtheorem{theorem}{Theorem}[section]

\newtheorem{claim}[theorem]{Claim}

\newtheorem{corollary}[theorem]{Corollary}

\newtheorem{definition}[theorem]{Definition}
\newtheorem{problem}[theorem]{Problem}

\newtheorem{conjecture}[theorem]{Conjecture}

\newcommand{\ignore}[1]{}

\begin{document}

\begin{spacing}{1.0}
\date{}
		\title{New bounds of two hypergraph Ramsey problems}
		
		\author{Chunchao Fan,\footnote{Center for Discrete Mathematics, Fuzhou University, Fuzhou, 350108, P.~R.~China. Email: {\tt 1807951575@qq.com}.}
            \;\; \;Xinyu Hu,\footnote{Data Science Institute, Shandong University, Jinan, 250100, P. R. China. Email: {\tt huxinyu@sdu.edu.cn}. }
			\;\; \; Qizhong Lin,\footnote{Center for Discrete Mathematics, Fuzhou University, Fuzhou, 350108, P.~R.~China. Email: {\tt linqizhong@fzu.edu.cn}. Supported in part  by National Key R\&D Program of China (Grant No. 2023YFA1010202), NSFC (No.\ 12171088, 12226401) and NSFFJ (No. 2022J02018).} \;\; \; Xin Lu\footnote{Zhengdong New District Yongfeng School, Zhengzhou, 451460, P.~R.~China. Email: {\tt lx20190902@163.com.}}
             }
\maketitle

\begin{abstract}

We focus on two hypergraph Ramsey problems.  First, we consider the Erd\H{o}s-Hajnal function $r_k(k+1,t;n)$.
In 1972, Erd\H{o}s and Hajnal conjectured that the tower growth rate of $r_k(k+1,t;n)$ is $t-1$ for each $2\le t\le k$. To finish this conjecture, it remains to show that the tower growth rate of $r_4(5,4;n)$ is three. We prove a superexponential lower bound for $r_4(5,4;n)$, which improves the previous best lower bound $r_4(5,4;n)\geq  2^{\Omega(n^2)}$ from Mubayi and Suk (\emph{J. Eur. Math. Soc., 2020}).  
Second, we prove an upper bound for the hypergraph Erd\H{o}s-Rogers function $f^{(k)}_{k+1,k+2}(N)$ that is an iterated $(k-3)$-fold logarithm in $N$ for each $k\geq 5$. This improves the previous upper bound that is an iterated $(k-13)$-fold logarithm in $N$ for $k\ge14$ due to Mubayi and Suk (\emph{J. London Math. Soc., 2018}), in which they conjectured that $f^{(k)}_{k+1,k+2}(N)$ is an iterated $(k-2)$-fold logarithm in $N$ for each $k\ge3$.

\medskip
\textbf{Keywords:} Hypergraph Ramsey number; Erd\H{o}s-Hajnal function; Erd\H{o}s-Rogers function

\end{abstract}

\section{Introduction}
We write $K^{(k)}_n$ for the complete $k$-uniform hypergraph ($k$-graph for short) on an $n$-element vertex set.
The Ramsey number $r_k(s,n)$ is the minimum integer $N$ such that for any red/blue edge coloring of the $K^{(k)}_N$, there is either a red copy of $K^{(k)}_s$ or a blue $K^{(k)}_n$. Estimating the Ramsey number $r_k(s,n)$ is a fundamental and notoriously difficult problem in combinatorics, which has attracted the most attention since 1935. 

For graphs, we know  $2^{n/2}<r_2(n,n)<2^{2n}$ from the two classical papers \cite{E-1,E-S}. We refer the reader to \cite{S-1,g-r,Thom,C-1,sah,C-G-M-S,G-N-N-L} for many improvements of the lower and upper bounds. In particular,  Conlon \cite{C-1} and subsequently Sah \cite{sah} made superpolynomial improvements. A recent breakthrough due to Campos, Griffiths,
Morris and Sahasrabudhe \cite{C-G-M-S} yields the first improvement to the exponential constant, which states that $r_2(n,n)<(4-\varepsilon)^n$ for some positive constant $\varepsilon$. More recently, Gupta, Ndiaye, Norin and Wei \cite{G-N-N-L} improved the upper bound to $3.8^{n+o(n)}$. 

However, for $3$-graphs, there still exists an exponential gap between the best known
upper and lower bounds for $r_3(n,n)$. Erd\H{o}s, Hajnal and Rado \cite{E-H-R,E-R-2} showed that 
$$2^{\Omega(n^2)}<r_3(n,n)<2^{2^{O(n)}}.$$
They conjectured that $r_3(n,n)>2^{2^{\Omega(n)}}$ and Erd\H{o}s offered a \$500 reward for a proof (see  \cite{Chung}). The upper bound has been improved to $2^{2^{(2+o(1))n}}$ by Conlon, Fox and Sudakov \cite{C-F-S}.

 For $k\ge4$, their results \cite{E-H-R,E-R-2} also yield an exponential gap between the lower and upper bounds for $r_k(n,n)$, i.e.,
$$\textrm{tw}_{k-1}(\Omega(n^2))<r_k(n,n)<\textrm{ tw}_k(O(n)),$$
where the tower function is defined recursively as $\twr_1(x)=x$ and $\twr_{i+1}(x)=2^{\twr_i(x)}$. 
Erd\H{o}s, Hajnal and Rado \cite{E-H-R} conjectured that  $r_k(n,n)=\twr_k(\Theta(n))$.

Off-diagonal Ramsey numbers $r_k(s,n)$ have also been extensively studied where $k$ and $s$ are fixed. For graphs, we know that $r_2(s,n)=n^{\Theta(1)}$ from \cite{sp75,S-1,A-K-S,she,lrz,B-1,B-K,kim,M-V}. In particular, we know  $r_2(3,n)=\Theta(n^2/\log n)$ from \cite{A-K-S,B-1,B-K,kim}, and $\Omega(n^3/\log ^4n)\le r_2(4,n)\le O(n^3/\log^2n)$ where the upper bound can be seen in \cite{A-K-S,lrz} while the lower bound was a recent breakthrough by Mattheus and Verstra\"{e}te \cite{M-V}. 

For $3$-graphs, improving the upper bound by Erd\H{o}s and Rado \cite{E-R-2} and the lower bound by Erd\H{o}s and Hajnal \cite{E-H-Con},  Conlon, Fox and Sudakov \cite{C-F-S} obtained that for $s\ge4$, $$2^{\Omega (n\log n)}\le r_3(s,n)\le2^{O(n^{s-2}\log n)}.$$

For $k\ge4$, Erd\H{o}s and Rado \cite{E-R-2} showed that $$r_k(s,n)\le \twr_{k-1}(n^{\Theta(1)}),$$ and Erd\H{o}s and Hajnal \cite{E-H-Con} conjectured that this bound is the correct tower growth rate for every fixed $k\geq 4$ and $s\geq k+1$. Conlon, Fox and Sudakov \cite{cfs} veriﬁed this conjecture for  $s\ge\lceil5k/2\rceil-3$. Subsequently, Mubayi and Suk \cite{M-S} and independently Conlon, Fox and Sudakov (unpublished) veriﬁed this conjecture for all $s\ge k+3$. 
Improving the bounds of $r_4(5,n)$ and $r_4(6,n)$ in \cite{M-S,M-S-2},  Mubayi and Suk \cite{M-S-1} verified the conjecture for $s\ge k+2$ by showing $r_4(6,n)\ge2^{2^{\Omega(n^{1/5})}}$ and consequently $r_k(k+2,n)\ge \twr_{k-1}(\Omega(n^{1/5}))$; and for the last case $s=k+1$, they proved $r_4(5,n)\ge2^{n^{\Omega(\log n)}},$ and consequently $r_k(k+1,n)\ge \twr_{k-2}(n^{\Omega(\log n)})$. To close the exponential gap between the best known lower and upper bounds for $r_k(k+1,n)$ for $k\ge4$, it suffices to show that $$r_4(5,n)\ge2^{2^{n^{\Theta(1)}}}$$ from the stepping-up lemma \cite[Lemma 2.1]{M-S-1}. Furthermore, such a lower bound would be true if we can show that $r_3(n,n)\ge\twr_3(\Omega(n))$ from \cite[Corollary 1.8]{M-S}. 

\medskip\noindent
{\bf Organization.}
We mainly consider two related hypergraph Ramsey problems: the hypergraph Erd\H{o}s-Hajnal problem and the hypergraph Erd\H{o}s-Rogers problem, which will be described in the relevant sections respectively. In Section \ref{B-p}, we will give some notations and basic properties we shall use throughout the paper. 
In Section \ref{e-h-p}, we consider the hypergraph Erd\H{o}s-Hajnal problem. Erd\H{o}s and Hajnal \cite{E-H-Con} conjectured that the tower growth rate of $r_k(k+1,t;n)$ is $t-1$ for each $2\le t\le k$. To complete this conjecture, it remains to show that the tower growth rate of $r_4(5,4;n)$ is three. We prove a superexponential lower bound for $r_4(5,4;n)$, which improves the previous best lower bound $r_4(5,4;n)\geq  2^{\Omega(n^2)}$ from Mubayi and Suk \cite[Theorem 6]{M-S-3} (as was pointed out in the nice survey \cite{M-S-4}). Consequently, we obtain new lower bounds for the hypergraph Erd\H{o}s-Hajnal function $r_k(k+1,k;n)$ for $k\ge4$, improving the previous best lower bound by Mubayi, Suk and Zhu \cite[Corollary 1.3]{M-S-Z}. 
In Section \ref{e-r-p}, we will focus on the hypergraph Erd\H{o}s-Rogers function $f^{(k)}_{k+1,k+2}(N)$. 
Conlon, Fox and Sudakov \cite{C-F-S-2} also proposed to close the gap between the upper and lower bounds of $f^{(k)}_{k+1,k+2}(N)$. Further, Mubayi and Suk \cite{M-S-2} conjectured that $f^{(k)}_{k+1,k+2}(N)$ is an iterated $(k-2)$-fold logarithm in $N$ for each $k\ge3$. We prove an upper bound for $f^{(k)}_{k+1,k+2}(N)$ that is an iterated $(k-3)$-fold logarithm in $N$ for each $k\geq 5$. This improves the previous upper bound that is an iterated $(k-13)$-fold logarithm in $N$ for $k\ge14$ due to Mubayi and Suk \cite[Theorem 3.2]{M-S-2}, and makes a substantial progress towards the conjecture of  Mubayi and Suk \cite{M-S-2}. To complete this conjecture, it only remains to show $g(4,2^N)=O\left((g(3,N))^{\Theta(1)}\right)$.

\section{Notations and basic properties}\label{B-p}

In this paper, we will apply several variants of the Erd\H os-Hajnal stepping-up lemma. Usually, it is not easy to construct suitable colorings when applying the stepping-up lemma. Building upon new constructions of the colorings, we are able to make progress of some hypergraph Ramsey problems.

We shall use the following notations and definitions  unless otherwise stated.
Let $\varphi$ be a red/blue coloring of the edges of the complete $(k-1)$-graph on the vertex set $\{0,1,\ldots, N-1\}$ satisfying some properties. We shall use $\varphi$ to define a red/blue coloring $\chi$ of the edges of the complete $k$-graph on the vertex set $V$, say $V=\{0,1,\ldots, 2^N-1\}$. We always use $(a_1,a_2,\ldots,a_{r})$ to denote a subset with $a_1<a_2<\cdots<a_{r}$.

For any $a \in V$, write $a=\sum_{i=0}^{N-1}a(i)2^i$ where $a(i) \in \{0,1\}$ for each $i$. For $a \not = b$, let $\delta(a,b)$ denote the largest $i$ such that $a(i) \not = b(i)$. Given any set $A=(a_1,a_2,\ldots,a_{r})$ of $V$, we always write for $1\le i\le r-1$, $$\delta_i=\delta(a_i,a_{i+1}),\;\;\text{and} \;\;\delta(A)=\{\delta_1,\delta_2,\dots,\delta_{r-1}\}.$$

We have the following stepping-up properties, see in \cite{grs}.

\begin{description}\label{p-ab}

\item[Property A:] For every triple $a < b < c$, $\delta(a,b) \not = \delta(b,c)$ .

\item[Property B:] For $a_1 < \cdots < a_r$, $\delta(a_1,a_{r}) = \max_{1 \leq i \leq r-1}\delta_i$.

\end{description}

We will also use the following stepping-up properties, which are easy consequences of Properties A and B.

\begin{description}\label{p-c}

\item[Property C:]  For $\delta(a_1,a_{r}) = \max_{1 \leq i \leq r-1}\delta_i$, there is a unique $\delta_i$ which achieves the maximum.
    
 \end{description} 
 
\noindent
{\bf Proof.} Suppose that there are indices $i<j$ such that
$\delta(a_i, a_{i+1})=\delta(a_{j}, a_{j+1})=\max_{1 \leq i \leq r-1} \delta_i.$
Then, by Property B, $\delta(a_i, a_{j})=\delta(a_{j}, a_{j+1})$, which is impossible from Property A.
\hfill$\Box$

\begin{description}\label{p-d}

\item[Property D:] For every 4-tuple $a_1 < \cdots < a_4$, if $\de1 > \de2$, then $\de1 \neq \de3$.   
    
 \end{description} 

\noindent
{\bf Proof.} Otherwise, suppose $\de1 = \de3$.
Then, by Property B, $\delta(a_1, a_{3})=\de1 =\de 3=\delta(a_{3}, a_{4})$. This contradicts Property A since $a_{1}<a_{3}<a_{4}$.
\hfill$\Box$

\medskip
Let $\varphi$ be defined as above. Suppose that we define $\chi$ as follows: for any $k$-tuple $(a_1,\ldots,a_k)$ of $V$, if $\delta_1,\ldots,\delta_{k-1}$ form a monotone sequence, i.e., $\delta_1<\ldots<\delta_{k-1}$ or $\delta_1>\ldots>\delta_{k-1}$, then  $$\chi(a_1,\ldots,a_{k})=\varphi(\delta_1,\ldots,\delta_{k-1}).$$

The following property is clear from Properties A and B.
\begin{description}\label{p-e}

\item[Property E:] For any $a_1 < \cdots < a_r$, suppose that $\delta_1,\ldots, \delta_{r-1}$ form a monotone sequence.  If $\chi$ colors every $k$-tuple in $(a_1,\ldots, a_r)$ red (blue), then $\varphi$ colors every $(k-1)$-tuple in $\{\delta_1,\ldots, \delta_{r-1}\}$ red (blue).

\end{description}

We can understand the red/blue edge coloring in another way: graph and complement graph. 
Let $H'$ be a $(k-1)$-graph satisfying some properties. Sometimes, it is convenient to apply $H'$ to define a $k$-graph $H$ on the vertex set $V$. 
Suppose we define the edges of $H$ as follows: for any $k$-tuple $(a_1,\ldots,a_k)$ of $V$, if $\delta_1,\ldots,\delta_{k-1}$ form a monotone sequence, then $(a_1,\ldots,a_k)\in E(H)$ if and only if $\{\delta_1,\ldots,\delta_{k-1}\}\in E(H')$. 

\medskip
We have the following property, which is a variant of Property E.

\begin{description}\label{p-f}

\item[Property F:]
For any $a_1 < \cdots < a_r$, suppose that $\delta_1,\ldots, \delta_{r-1}$ form a monotone sequence.  If every $k$-tuple in $(a_1,\ldots, a_r)$ is in $E(H)$ (in $\overline{E}(H)$), then every $(k-1)$-tuple in $\{\delta_1,\ldots, \delta_{r-1}\}$ is in $E(H')$ (in $\overline{E}(H')$).

\end{description}

For any set $A = (a_1,\ldots, a_r)$, we call $\delta_i$ a {\it local minimum} of $A$ if $\delta_{i-1}>\delta_i<\delta_{i+1}$, a {\it local maximum} if $\delta_{i-1}<\delta_i>\delta_{i+1}$, and a {\it local extremum} if it is either a local minimum or a local maximum.  We call $\delta_i$ a \emph{local monotone} if $\delta_{i-1} < \delta_i < \delta_{i + 1}$ or $\delta_{i-1} > \delta_i > \delta_{i + 1}$.
 
 For any set $A = (a_1,\ldots, a_r)$ where $r\ge4$, 
  let $m(A)$ and $n(A)$ denote the number of {\em local extrema} and {\em local monotone} in the sequence $\delta_1,\ldots, \delta_{r-1}$, respectively. Since for $2\le i\le r-2$, each $\delta_i$ is either a local extremum or a local monotone, we have $m(A)+n(A)=r-3.$

\begin{description}\label{p-g}

\item[Property G:]
Any two consecutive local maxima are distinct.

\end{description}

\noindent
{\bf Proof.} Let $\de p$ and $\de q$  where $p<q$ be any two consecutive local maxima. 
On the contrary, suppose $\de p = \de q$.
Then, by Property B, $\delta(a_p, a_{q})=\de p =\de q=\delta(a_{q}, a_{q+1})$. This contradicts Property A since $a_{p}<a_{q}<a_{q+1}$.
\hfill$\Box$

\medskip
For any set $A=(a_1,\dots,a_r)$,  we always write $A\setminus a_i$ for the subset of $A$ by deleting $a_i$, and write $A\setminus (a_i,a_j)$ for the subset of $A$ by deleting the elements of $a_i$ and $a_j$, respectively. Write 
$[A]_i^j=(a_i,\dots,a_j).$ In the following, let all logarithms be base 2.

\section{Hypergraph Erd\H{o}s-Hajnal problem}\label{e-h-p}
Since the Ramsey number $r_k(k+1,n)$ for $k\ge4$ is not very well understood, Erd\H{o}s and Hajnal \cite{E-H-Con} introduced the following
more general function.

\begin{definition}
For integers $2\le k<n$ and $2\le t\le k+1$, let $r_k(k+1,t;n)$ be the minimum $N$
such that for every red/blue edge coloring of the $K^{(k)}_N$, there is a $k$-graph on $k+1$ vertices with at least $t$ edges colored red, or there exists a blue copy of $K^{(k)}_n$.
\end{definition}

It is clear that $r_k(k+1,1;n)=n$, and $r_k(k+1,k+1;n)=r_k(k+1,n)$. For each $2\le t \le k$, Erd\H{o}s and
Hajnal \cite{E-H-Con} showed that $r_k(k+1,t;n)<\twr_{t-1}(n^{\Theta(1)}),$ and further Erd\H{o}s and Hajnal conjectured that the upper bound is the correct tower growth rate for every fixed $k\geq 4$ and $2\le t\le k$. Namely,
\begin{align}\label{E-H-conj}
 r_k(k+1,t;n)=\twr_{t-1}(n^{\Theta(1)}).
\end{align}
In particular, they conjectured that $r_k(k+1,k;n)$ has the same tower growth rate as $r_k(k+1,n)$.

We know that Erd\H{o}s and Hajnal's conjecture (\ref{E-H-conj}) was verified to be true for all cases but the last case where $k\ge4$ and $t=k$:

 \medskip
$\bullet$ for $k\le3$ or $t\le 3$, see Erd\H{o}s and Hajnal \cite{E-H-Con}.

 \medskip
$\bullet$ for $k\ge5$ and $3\le t\le k-2$, see Mubayi and Suk \cite{M-S-3} in a stronger form.

 \medskip
$\bullet$ for $k\ge4$ and $t=k-1$, see Mubayi, Suk and  Zhu \cite{M-S-Z}.

 \medskip

It remains to verify the last case of the conjecture.

 \begin{conjecture}[Erd\H{o}s and Hajnal \cite{E-H-Con}]\label{conj}
For $k\ge4$,
$$r_k(k+1,k;n)=\twr_{k-1}(n^{\Theta(1)}).$$
\end{conjecture}

 For the general lower bound of $r_k(k+1,t;n)$,  Mubayi and Suk \cite{M-S-3} proved the following iteration result. In particular, the case where $t=k+1$ was also proved by Mubayi and Suk \cite{M-S-1}.
\begin{theorem}[Mubayi and Suk \cite{M-S-1,M-S-3}]\label{int}
For $k\ge 6$ and $5\le t\le k+1$, $$r_k(k+1,t;2kn)>2^{r_{k-1}(k,t-1;n)-1}.$$
\end{theorem}


From Theorem \ref{int} together with that $r_5(6,5;4n^2)>2^{r_4(5,4;n)-1}$ by Mubayi, Suk and Zhu \cite{M-S-Z}, we are able to confirm Conjecture \ref{conj} if we can prove $r_4(5,4;n)\geq  2^{2^{n^{\Theta(1)}}}$. It would be a challenge since this already implies that $r_4(5,n)\geq  2^{2^{n^{\Theta(1)}}}$.

The current best lower and upper bounds of $r_4(5,4;n)$ are as follows:
\begin{align}\label{o-r454}
2^{\Omega(n^2)}<r_4(5,4;n)<2^{2^{n^{\Theta(1)}}},
\end{align}
where the lower bound follows from a result of Mubayi and Suk \cite[Theorem 6]{M-S-3} as was pointed out in the nice survey \cite{M-S-4} by noting that $r_4(5,4;n)\ge r_4(5,3;n)=2^{\Omega(n^2)}$,  while the upper bound was already known from Erd\H{o}s and
Hajnal \cite{E-H-Con}.

In general, Mubayi, Suk and Zhu \cite[Corollary 1.3]{M-S-Z} proved that for each $k\ge5$, 
\begin{align}\label{o-g-l}
r_k(k+1,k;n)>\twr_{k-2}(n^{\Theta(1)}).
\end{align}


In this paper, we obtain a superexponential lower bound for $r_4(5,4;n)$,
which improves the previous best lower bound in (\ref{o-r454}) due to Mubayi and Suk \cite{M-S-3} as follows.

\begin{theorem}\label{5,4;n}
For all sufficiently large $n$,
\[
r_4(5,4;n) \geq  2^{n^{c \log\log n}},
\]where $c$ is a positive constant.
\end{theorem}

Generally, together with the fact that $r_5(6,5;4n^2)>2^{r_4(5,4;n)-1}$ and Theorem \ref{int}, we have the following result, which improves (\ref{o-g-l}) by Mubayi, Suk and Zhu \cite{M-S-Z}.

\begin{corollary}\label{k+1,k;n}
For each $k\ge 5$, $r_k(k+1,k;n)>\twr_{k-2}(n^{\Omega(\log\log n)}).$
\end{corollary}

\subsection{Proof of Theorem \ref{5,4;n}}


From Erd\H{o}s and Hajnal \cite{E-H-Con}, we know $2^{\Omega(n)}\le r_3(4,3;n)\le 2^{O(n\log n)}.$
As the smallest 3-graph whose Ramsey number with a clique is at least exponential, it has attracted much of attention, see \cite{C-F-S,E-2,M-S-2,M-S-4}.
Fox and  He  \cite[Theorem 1.1]{F-H} determined the order of the exponent.
\begin{theorem}[Fox and  He \cite{F-H}]\label{4,3;n}
For all sufficiently large $n$, $r_3(4,3;n)=2^{\Theta(n\log n)}.$
\end{theorem}


Now, Theorem \ref{5,4;n} follows from the following result immediately.
\begin{theorem}\label{lm-5,4;n}
For all sufficiently large $n$, $r_4(5,4;2^{2n})> 2^{r_3(4, 3; n) - 1}$.
\end{theorem}

\noindent\textbf{Proof.} We will apply a variant of the Erd\H os-Hajnal stepping-up lemma. 
 Let $\varphi$ be a red/blue coloring of the edges of the complete $3$-graph on the vertex set $\{0,1,\ldots, r_3\left(4, 3; n\right) - 2\}$ such that there are at most $2$ red edges
among every $4$ vertices and there is no blue $K_{n}^{(3)}$. Let $N= 2^{r_3(4,3; n) - 1}$. We shall use $\varphi$ to define a red/blue coloring $\chi$ of the edges of $K_N^{(4)}$ on the vertex set $V=\{0,1,\ldots, N-1\}$.
For any 4-tuple $(a_1, a_2, a_3, a_4)\subseteq V$,  we define $\chi$ as follows:

\begin{equation}
 \chi(a_1, a_2, a_3, a_4) =\begin{cases}

\varphi(\delta_1,\delta_2,\delta_{3}) \qquad \hbox{ if $\delta_1,\delta_2,\delta_{3}$ are monotone,} \cr
blue \qquad \qquad \quad \hbox{\, otherwise.} \notag
\end{cases}
\end{equation}
Therefore, Properties A, B and C hold, and Property E holds where $k=4$ (see in Section \ref{B-p}).

In the following, we will show that there are at most 3 red edges among every $5$ vertices and no blue $K_{m}^{(4)}$ in coloring $\chi$ where $m=2^{2n}$.

Suppose first that there exist $5$ vertices $a_1,\dots,a_5$ where $a_1 < a_2 < \cdots < a_5$ that contain at least $4$ red edges. We shall show that this will lead to a contradiction.  Let $e_i = (a_1,\dots,a_5)\setminus a_i$. Let $\dd(e_i)$ be the resulting sequence of $\dd$'s. In particular, for $i = 1$, $\dd(e_1) = (\de2,\de3,\de4)$. For $2 \leq i \leq 4$, $\dd(e_i) = (\de1, \dots, \dd(a_{i-1},a_{i+1}), \dots, \de4)$. For $i = 5$, $\dd(e_5) = (\de1,\de2,\de3)$. For clarity, we will write $\dd(e_1) = (\de2>\de3,\de4)$ if we also know $\de2>\de3$.
 In the following, we will often use the fact that if $2 \leq i \leq 4$, then $\delta(a_{i-1},a_{i+1}) = \max\{\de{i-1},\de{i}\}$ by Property B.

Now we will consider cases depending on the ordering of $\de1,\de2,\de3,\de4$, and there are 8 possible orderings which will be split into five cases  as follows.

\medskip

\textit{Case 1}. $\de1 > \de2 > \de3 > \de4$. 
This implies that
\begin{align*}
\dd(e_1) &= (\de2 > \de3 > \de4 ),\\
\dd(e_2) &= (\de1 >\de3  > \de4 ),\\
\dd(e_3) &= (\de1 > \de2 > \de4),\\
\dd(e_4) =\dd(e_5)&= (\de1 > \de2 > \de3).
\end{align*}
Since there are at least $4$ red edges among these $5$ vertices, $e_4,e_5$ are red and at most one of $e_1,e_2,e_3$ is blue, say $e_1$, without loss of generality. This implies that $\varphi(\de1, \de3, \de4)$, $\varphi(\de1, \de2, \de4)$ and $\varphi(\de1, \de2, \de3)$ are red among these $4$ vertices, a contradiction.

\medskip

\textit{Case 2}. $\de1 < \de2 < \de3 < \de4$. 
This implies that
\begin{align*}
\dd(e_1) = \dd(e_2) &= (\de2 < \de3 < \de4 ),\\
\dd(e_3) &= (\de1 < \de3 < \de4),\\
\dd(e_4) &= (\de1 < \de2 <\de4),\\
\dd(e_5) &= (\de1 < \de2 < \de3).
\end{align*}
Similarly, since there are at least $4$ red edges among these $5$ vertices, $e_1$ and $e_2$ must be red  and at most one of $e_3,e_4,e_5$ is blue, say $e_5$, without loss of generality. This implies that $\varphi(\de2, \de3, \de4)$, $\varphi(\de1, \de3, \de4)$ and $\varphi(\de1, \de2, \de4)$ are red among these $4$ vertices, a contradiction.

\medskip

\textit{Case 3}. $\de1 > \de2 > \de3 < \de4 $, or $\de1 < \de2 < \de3 > \de4 $. 
In the first case, $\dd(e_1)=(\de2 >\de3 < \de4 )$ and  $\dd(e_2)=(\de 1 > \de3 < \de4)$, so both $e_1$ and $e_2$ are blue, contradicting the fact that there are at least 4 red edges among these 5 vertices. In the second case, $\dd(e_1)=\dd(e_2)=(\de 2 < \de3 > \de4)$, so both $e_1$ and $e_2$ are blue, again leading to a contradiction.

\medskip

\textit{Case 4}. $\de1 > \de2 < \de3 > \de4 $, or $\de1 < \de2 > \de3 < \de4 $. 
In the first case, we obtain that $$\dd(e_4)=\dd(e_5)=(\de 1 > \de2 < \de3),$$ so $e_4$ and $e_5$ are blue, a contradiction. In the second case, $\dd(e_1)=\dd(e_2)=(\de 2 > \de3 < \de4)$, so $e_1$ and $e_2$ are blue, again leading to a contradiction.

\medskip

\textit{Case 5}. $\de1 > \de2 < \de3 < \de4 $, or $\de1 < \de2 > \de3 > \de4 $.
 In the first case, $\dd(e_4)=(\de 1 > \de2 < \de 4)$, and $\dd(e_5)=(\de 1 > \de2 < \de 3)$, so both $e_4$ and $e_5$ are blue, a contradiction. In the second case, $\dd(e_4)=\dd(e_5)=(\de 1 < \de2 > \de3)$, we can also deduce that both $e_4$ and $e_5$ are blue, again leading to a contradiction.


\begin{figure}
\begin{center}

\tikzset{every picture/.style={line width=0.70pt}} 

\begin{tikzpicture}[x=0.70pt,y=0.70pt,yscale=-1,xscale=1]

\draw    (95,154.88) -- (572,154.88) (106,150.88) -- (106,158.88)(139,150.88) -- (139,158.88)(172,150.88) -- (172,158.88)(205,150.88) -- (205,158.88)(238,150.88) -- (238,158.88)(271,150.88) -- (271,158.88)(304,150.88) -- (304,158.88)(337,150.88) -- (337,158.88)(370,150.88) -- (370,158.88)(403,150.88) -- (403,158.88)(436,150.88) -- (436,158.88)(469,150.88) -- (469,158.88)(502,150.88) -- (502,158.88)(535,150.88) -- (535,158.88)(568,150.88) -- (568,158.88) ;
\draw [color={rgb, 255:red, 208; green, 2; blue, 27 }  ,draw opacity=1 ]   (187.78,148.71) -- (187.75,129.96) ;
\draw [shift={(187.75,127.96)}, rotate = 89.89] [color={rgb, 255:red, 208; green, 2; blue, 27 }  ,draw opacity=1 ][line width=0.75]    (10.93,-3.29) .. controls (6.95,-1.4) and (3.31,-0.3) .. (0,0) .. controls (3.31,0.3) and (6.95,1.4) .. (10.93,3.29)   ;
\draw [color={rgb, 255:red, 208; green, 2; blue, 27 }  ,draw opacity=1 ]   (287.78,146.71) -- (287.75,127.96) ;
\draw [shift={(287.75,125.96)}, rotate = 89.89] [color={rgb, 255:red, 208; green, 2; blue, 27 }  ,draw opacity=1 ][line width=0.75]    (10.93,-3.29) .. controls (6.95,-1.4) and (3.31,-0.3) .. (0,0) .. controls (3.31,0.3) and (6.95,1.4) .. (10.93,3.29)   ;
\draw [color={rgb, 255:red, 74; green, 144; blue, 226 }  ,draw opacity=1 ]   (519.78,147.71) -- (519.75,128.96) ;
\draw [shift={(519.75,126.96)}, rotate = 89.89] [color={rgb, 255:red, 74; green, 144; blue, 226 }  ,draw opacity=1 ][line width=0.75]    (10.93,-3.29) .. controls (6.95,-1.4) and (3.31,-0.3) .. (0,0) .. controls (3.31,0.3) and (6.95,1.4) .. (10.93,3.29)   ;
\draw [color={rgb, 255:red, 74; green, 144; blue, 226 }  ,draw opacity=1 ]   (451.78,147.71) -- (451.75,128.96) ;
\draw [shift={(451.75,126.96)}, rotate = 89.89] [color={rgb, 255:red, 74; green, 144; blue, 226 }  ,draw opacity=1 ][line width=0.75]    (10.93,-3.29) .. controls (6.95,-1.4) and (3.31,-0.3) .. (0,0) .. controls (3.31,0.3) and (6.95,1.4) .. (10.93,3.29)   ;
\draw  [color={rgb, 255:red, 208; green, 2; blue, 27 }  ,draw opacity=1 ] (206.5,198.5) .. controls (206.49,203.17) and (208.81,205.51) .. (213.48,205.52) -- (378.52,205.98) .. controls (385.19,206) and (388.51,208.34) .. (388.5,213.01) .. controls (388.51,208.34) and (391.85,206.02) .. (398.52,206.04)(395.52,206.03) -- (561.98,206.49) .. controls (566.65,206.5) and (568.99,204.18) .. (569,199.51) ;
\draw  [color={rgb, 255:red, 74; green, 144; blue, 226 }  ,draw opacity=1 ] (568.5,81.5) .. controls (568.52,76.83) and (566.2,74.49) .. (561.53,74.47) -- (449.46,74.04) .. controls (442.79,74.01) and (439.47,71.67) .. (439.49,67) .. controls (439.47,71.67) and (436.13,73.99) .. (429.46,73.97)(432.46,73.98) -- (310.03,73.51) .. controls (305.36,73.49) and (303.02,75.81) .. (303,80.48) ;
\draw    (352.78,146.71) -- (352.75,127.96) ;
\draw [shift={(352.75,125.96)}, rotate = 89.89] [color={rgb, 255:red, 0; green, 0; blue, 0 }  ][line width=0.75]    (10.93,-3.29) .. controls (6.95,-1.4) and (3.31,-0.3) .. (0,0) .. controls (3.31,0.3) and (6.95,1.4) .. (10.93,3.29)   ;
\draw   (304.5,176.5) .. controls (304.47,181.17) and (306.78,183.52) .. (311.45,183.55) -- (359.54,183.92) .. controls (366.21,183.97) and (369.52,186.33) .. (369.49,191) .. controls (369.52,186.33) and (372.87,184.03) .. (379.54,184.08)(376.54,184.05) -- (428.95,184.46) .. controls (433.62,184.49) and (435.97,182.18) .. (436,177.51) ;
\draw  [color={rgb, 255:red, 208; green, 2; blue, 27 }  ,draw opacity=1 ] (105.5,222.5) .. controls (105.51,227.17) and (107.84,229.5) .. (112.51,229.49) -- (329.41,229.03) .. controls (336.08,229.02) and (339.41,231.34) .. (339.42,236.01) .. controls (339.41,231.34) and (342.74,229) .. (349.41,228.99)(346.41,229) -- (564.01,228.54) .. controls (568.68,228.53) and (571.01,226.2) .. (571,221.53) ;
\draw  [color={rgb, 255:red, 74; green, 144; blue, 226 }  ,draw opacity=1 ][line width=0.75]  (501,102.5) .. controls (501,97.83) and (498.67,95.5) .. (494,95.5) -- (413.5,95.5) .. controls (406.83,95.5) and (403.5,93.17) .. (403.5,88.5) .. controls (403.5,93.17) and (400.17,95.5) .. (393.5,95.5)(396.5,95.5) -- (311.5,95.5) .. controls (306.83,95.5) and (304.5,97.83) .. (304.5,102.5) ;

\draw (561.45,160) node [anchor=north west][inner sep=0.75pt]   [align=left] {$\displaystyle \mathnormal{a}$$\displaystyle \mathnormal{_{m}}$\textit{\textsubscript{}}};
\draw (195.55,160) node [anchor=north west][inner sep=0.75pt]   [align=left] {$\displaystyle \mathnormal{a}$$\displaystyle _\mathnormal{s_{1} +1}$\textit{\textsubscript{}}};
\draw (166.08,160) node [anchor=north west][inner sep=0.75pt]   [align=left] {$\displaystyle \mathnormal{a}$$\displaystyle \mathnormal{_{s_{1}}}$\textit{\textsubscript{}}};
\draw (99.28,160) node [anchor=north west][inner sep=0.75pt]   [align=left] {$\displaystyle \mathnormal{a}$$\displaystyle \mathnormal{_{1}}$\textit{\textsubscript{}}};
\draw (181.24,107.19) node [anchor=north west][inner sep=0.75pt]  [font=\small,color={rgb, 255:red, 208; green, 2; blue, 27 }  ,opacity=1 ,rotate=-359.06] [align=left] {$\displaystyle {\delta' _{1}}$\textit{\textsubscript{}}};
\draw (280.24,106.19) node [anchor=north west][inner sep=0.75pt]  [font=\small,color={rgb, 255:red, 208; green, 2; blue, 27 }  ,opacity=1 ,rotate=-359.06] [align=left] {$\displaystyle {\delta' _{2}}$\textit{\textsubscript{}}};
\draw (510.24,106.19) node [anchor=north west][inner sep=0.75pt]  [font=\small,color={rgb, 255:red, 74; green, 144; blue, 226 }  ,opacity=1 ,rotate=-359.06] [align=left] {$\displaystyle {\delta' _{3}}$\textit{\textsubscript{}}};
\draw (443.24,107.19) node [anchor=north west][inner sep=0.75pt]  [font=\small,color={rgb, 255:red, 74; green, 144; blue, 226 }  ,opacity=1 ,rotate=-359.06] [align=left] {$\displaystyle {\delta' _{4}}$\textit{\textsubscript{}}};
\draw (343.24,104.19) node [anchor=north west][inner sep=0.75pt]  [font=\small,rotate=-359.06] [align=left] {$\displaystyle {\delta' _{2n}}$\textit{\textsubscript{}}};
\draw (330.08,160) node [anchor=north west][inner sep=0.75pt]   [align=left] {$\displaystyle \mathnormal{a}$$\displaystyle \mathnormal{_{s_{2n}}}$\textit{\textsubscript{}}};
\draw (362.5,160) node [anchor=north west][inner sep=0.75pt]   [align=left] {$\displaystyle \mathnormal{a}$$\displaystyle \mathnormal{_{s_{2n} +1}}$\textit{\textsubscript{}}};
\draw (525,160) node [anchor=north west][inner sep=0.75pt]   [align=left] {$\displaystyle \mathnormal{a}$$\displaystyle \mathnormal{_{m-1}}$\textit{\textsubscript{}}};
\draw (307,110) node [anchor=north west][inner sep=0.75pt]   [align=left] {......};
\draw (390,110) node [anchor=north west][inner sep=0.75pt]   [align=left] {......};

\end{tikzpicture}

\end{center}
\begin{center}
Fig. 1. The search process.
\end{center}
\end{figure}

\medskip
Next we show that there is no blue $K_{m}^{(4)}$ in coloring $\chi$ where $m=2^{2n}$. 
On the contrary, suppose that there are $m$-tuple $A=(a_1,\ldots, a_{m})$ such that $\chi$ colors every $4$-tuple of $A$ blue. Recall $[A]_{i}^j=(a_i,\dots, a_{j})$, and $\delta([A]_{i}^j)=\{\delta_i,\dots,\delta_{j-1}\}$ for $1\le i<j\le m$.   
Our aim is to find a monotone sequence $\delta'$'s of size $t$ according to a sequence $(a_{i_1},\ldots, a_{i_{t+1}})$ such that $\delta'_l=\delta(a_{i_l},a_{i_{l+1}})$. The search process is as follows. We first find the unique largest $\delta_i$ in $\delta(A)$ as $\delta'_1=\delta(a_{s_1},a_{s_{1}+1})$. If $|[A]_{1}^{s_1}|\geq |[A]_{s_1+1}^m|$, then we find the unique largest $\delta_i$ in $\delta([A]_{1}^{s_1})$ as $\delta'_2=\delta(a_{s_2},a_{s_{2}+1})$; otherwise we find the unique largest $\delta_i$ in $\delta([A]_{s_1+1}^m)$ as $\delta'_2=\delta(a_{s_2},a_{s_{2}+1})$. 
The uniqueness of $\delta'_i$ follows from Property C.
We always search the next $\delta'_{l+1}$ in the larger set of $\{\delta([A]_{i}^{s_l}),\delta([A]_{s_l+1}^j)\}$  resulting from  $\delta'_{l}$. 
The number of vertices of each step for searching $\delta'_{l+1}$ is at most halved.
Therefore, the process can not stop until we find $\delta'_{2n}$. The search process can be seen in Figure 1.

By the pigeonhole principle, at least $n$ elements of $\{\delta'_1,\ldots,\delta'_{2n}\}$ are on the same side as $a_{s_{2n}}$. Without loss of generality, there exists a monotone decreasing sequence $\{\delta''_1,\ldots,\delta''_n\}$ where $\delta''_{l}=\delta(a'_{s_l},a'_{s_{l}+1})$ on the left of $a_{s_{2n}}$. By Property B, we find a sequence $(a'_{s_1},\ldots, a'_{s_{n}}, a'_{s_{n}+1})$ such that $\delta''_l=\delta(a'_{s_l},a'_{s_{l}+1})=\delta(a'_{s_l},a'_{s_{l+1}})$ for $1\leq l\leq n-1$ and $\delta''_n=\delta(a'_{s_n},a'_{s_{n}+1})$.
By Property E with $k=4$, $\varphi$ colors every 3-tuple from this set $\{\delta''_1,\ldots,\delta''_n\}$ blue, which is a contradiction.

Therefore, there are at most 3 red edges among every $5$ vertices and no blue $K_{m}^{(4)}$ in coloring $\chi$. The lower bound follows as desired. \hfill$\Box$

\section{Hypergraph Erd\H{o}s-Rogers problem}\label{e-r-p}

In this section, we focus on the hypergraph Erd\H{o}s-Rogers problem.
We first introduce the generalization of independent sets. For any integer $s \ge 2$, an $s$-independent set in a $k$-graph $H$ is a vertex subset that is $K^{(k)}_s$-free. So if
$s=k$, then it is just an independent set. Let $\alpha_s(H)$ denote the size of the largest $s$-independent set in $H$. 
Erd\H{o}s and Rogers \cite{E-R-1} proposed to determine the following function.
\begin{definition}
For $k\le s<t<N$, the function $f^{(k)}_{s,t}(N)$ is the minimum of $\alpha_s(H)$ taken over all $K^{(k)}_t$-free $k$-graphs $H$ of order $N$.
\end{definition}

We have the following:

\medskip
$\bullet$ $f^{(k)}_{s,t}(N)\ge n$ $\Leftrightarrow$ Every $K^{(k)}_t$-free $k$-graph $H$ of order $N$ contains a set of $n$ vertices that is $K^{(k)}_s$-free, i.e., $\alpha_s(H)\ge n$.

\medskip
$\bullet$ $f^{(k)}_{s,t}(N)<n$ $\Leftrightarrow$ There exists a $K^{(k)}_t$-free $k$-graph $H$ of order $N$ satisfying that each set of $n$ vertices contains a copy of $K^{(k)}_s$, i.e., $\alpha_s(H)<n$.

\medskip
Recall that the Ramsey number $r_k(t,n)$ is the minimum $N$ such that every red/blue coloring of the edges of $K^{(k)}_N$ contains a red $K^{(k)}_t$ or a blue $K^{(k)}_n$. We know that the Erd\H{o}s-Rogers function and the Ramsey number have the following relation:
\begin{align*}
r_k(t,n)=\min\left\{N:f^{(k)}_{k,t}(N)\ge n\right\}.
\end{align*}

Erd\H{o}s and Rogers \cite{E-R-1} initiated to study the function $f^{(k)}_{s,t}(N)$ for $k=2$ and $t=s+1$,
which was subsequently addressed by Bollob\'{a}s and Hind \cite{bh}, Krivelevich \cite{kri-1,kri-2}, Alon and Krivelevich \cite{ak}, Dudek and R\"{o}dl \cite{dr}, Wolfovitz \cite{W-1} and Dudek, Retter and R\"{o}dl \cite{D-R-R}. The previous best lower and upper bounds are as follows:
for each $s\ge3$,
\begin{align}\label{e-f-r}
\Omega\left(N^{1/2}\left(\frac{\log N}{\log \log N}\right)^{1/2}\right)<f^{(2)}_{s,s+1}(N)<O\left(N^{1/2}\log N\right),
\end{align}
where the lower bound is due to Dudek and Mubayi \cite{D-M} while the upper bound comes from Mubayi and Verstra\"{e}te \cite{M-V-1}.
We refer the reader to Mubayi and Suk \cite{M-S-4} for a nice survey on this topic, and \cite{GJ20,Janzer,H-L,Sud-1,Sud-2} etc. for some related results.

For $k\ge3$, it is much more difficult to evaluate the Erd\H{o}s–Rogers function. In the following, we will focus on the case where $s=k+1$.
For simplicity, we write $$g(k,N)=f^{(k)}_{k+1,k+2}(N).$$
Thus (\ref{e-f-r}) can be rephrased as $g(2,N)=N^{1/2+o(1)}$.

 For $k\ge3$, Dudek and Mubayi \cite{D-M} proved that
\begin{align}\label{mb-up-8}
\Omega\left((\log_{(k-2)} N)^{1/4}\right)<g(k,N)<O\left((\log N)^{1/(k-2)}\right),
\end{align}
 where $\log_{(i)}$ is the $\log$ function iterated $i$ times. 
 
 Conlon, Fox and Sudakov \cite{C-F-S-2} improved the lower bound to that
  \[
g(k,N)>\Omega\left((\log_{(k-2)} N)^{1/3-o(1)}\right),
  \]
  and they also proposed to close the gap between the upper and lower bounds. This was firstly achieved by Mubayi and Suk \cite[Theorem 3.2]{M-S-2}.
 
 \begin{theorem}[Mubayi and Suk \cite{M-S-2}]\label{mb-up}
For each $k\ge 14$,
$
g(k,N)=O(\log_{(k-13)}N).
$

\end{theorem}
 
Furthermore, they proposed to determine the correct number of iterations (which may well be $k-2$).

\begin{conjecture}[Mubayi and Suk \cite{M-S-2}]\label{con}
For each $k\ge 3$, $g(k,N)=(\log_{(k-2)}N)^{\Theta(1)}.$
\end{conjecture}

In this paper, building upon several new constructions of the colorings, we are able to obtain new upper bounds of $g(k,N)$ for all $k\ge5$ as follows.
\begin{theorem}\label{cr-main-2}
For each fixed $k \ge 5$, 
 $g(k,N) = O( \log_{(k-3)} N).$
\end{theorem}

{\em Remark.} Our result improves Theorem \ref{mb-up} and the upper bound in (\ref{mb-up-8}) for $k\ge5$. In particular, this makes a substantial progress towards Conjecture \ref{con}. In order to complete Conjecture \ref{con}, it only remains to show that $$g(4,2^N)=O\left((g(3,N))^{\Theta(1)}\right).$$

In order to prove Theorem \ref{cr-main-2}, it suffices to prove the following two iterations.

\begin{theorem}\label{Erdos-Roger}
We have $g(5,2^N)<4g^2(4,N)$. 
\end{theorem}

For each $k\ge6$, we have the following iteration.
\begin{theorem}\label{main-2}
For each fixed $k \ge 6$, $g(k,2^N)< (k+2) g(k-1, N)$. 
\end{theorem}

We first prove Theorem \ref{cr-main-2}, which follows from the above two theorems immediately.

\medskip\noindent
{\bf Proof of Theorem \ref{cr-main-2}.} From Theorem \ref{Erdos-Roger} and the upper bound in (\ref{mb-up-8}) for $k=4$, we know that $g(5,N)\le g(5, 2^{\lceil \log N \rceil} )=O(\log\log N).$ Together with Theorem \ref{main-2} we can inductively obtain that for $k\ge6$,
\begin{align*}
g(k,N) \le g(k, 2^{\lceil \log N \rceil} )
&< (k+2) \cdot g(k-1, \lceil \log N \rceil) 
\\&\le O( \log_{(k-4)} \lceil\log N \rceil)
\\&\leq O( \log_{(k-3)} N).
\end{align*}
The assertion follows.
\hfill$\Box$

\medskip
In the following, we shall give proofs for Theorem \ref{Erdos-Roger} and Theorem \ref{main-2}.
We mainly apply several variants of the Erd\H{o}s–Hajnal stepping-up lemma by building new constructions of the colorings.

\subsection{Proof of Theorem \ref{Erdos-Roger}}

We will show that $g(5,2^N)<4g^2(4,N)$. 
 From the definition of $g(4,N)$, there is a 4-graph $H'$ on the vertex set $\{0,1,\dots,N-1\}$ with $K^{(4)}_6$-free and $\alpha_5(H')=g(4,N)$. We shall use $H'$ to produce a $K^{(5)}_7$-free 5-graph $H$ on $V(H)=\{0,1,\dots,2^N-1\}$ with $\alpha_6(H)<4(\alpha_5(H'))^2$. We define the edges of $H$ as follows.  For any $5$-tuple $e=(a_1,\ldots,a_5)$ of $V(H)$, set $e\in E(H)$ if and only if one of the following holds:
\medskip

(I) $m(e)=0$, i.e., $\delta_1,\ldots, \delta_4$ form a monotone sequence and $\{\delta_1,\ldots, \delta_4\}\in E(H')$;
\medskip

(II) $m(e)=2$.



\medskip

Therefore, Properties A, B, D, F and G hold (see in Section 2).

We first show that $H$ contains no $7$-clique. On the contrary, there is a set $A=(a_1,\ldots,a_7)$ that induces a $K^{(5)}_7$ in $H$. 
For any $f\subseteq V(H)$, recall $m(f)$ and $n(f)$ are the number of {\em local extrema} and {\em local monotone} in $\delta(f)$, respectively.
\begin{claim} \label{mono} 
Let $f = (a_1,\ldots,  a_6)$, if $m(f) = 3$ and $\de 2$ is a local maximum, then $n(f\setminus a_i)=1$ for each $3\leq i\leq 4$.
In particular, either $\de 2$ or $\de 4$  is a local monotone.
\end{claim}

\noindent
{\bf Proof.} Since $\delta_2$ is a local maximum, we must have $\delta_1< \delta_2 > \delta_3 < \delta_4 >\delta_5$. Therefore, 
$$\dd(f\setminus a_3) =\dd(f\setminus a_4) =(\de1 < \de2,\ \de4>\de 5 ).$$ 
Since $\de2>\de3$, by Property D we have $\de2\neq\de 4$.
If $\de2<\de 4$, then $\de 2$ is a local monotone.
If $\de2>\de 4$, then $\de 4$ is a local monotone.
\hfill$\Box$

\medskip

\emph{Case 1.} $m(A)=0$. Then $\delta_1,\ldots, \delta_6 $ form a monotone sequence, and thus every 4-tuple in the set $\delta(A)=\{\delta_1,\ldots, \delta_6\}$ is an edge in $H'$ by Property F, contradicting the fact that $H'$ is $K^{(4)}_6$-free.

\medskip

\emph{Case 2.} $m(A)\ge 1$. Let $\delta_p$ be the first local extremum in the sequence $\delta(A)$, where $2\le p\le5$. If $3\le p\le5$, then the 5-tuple $e=(a_{p-2},a_{p-1},a_p,a_{p+1},a_{p+2})$ satisfies that $m(e)=1$ since $\delta_p$ is the first local extremum in the sequence $\delta(A)$. Thus $e\notin H$, a contradiction.
Therefore,  we may assume that $\delta_2$ is the first local extremum in the sequence $\delta(A)$.

Suppose $\delta_2$ is a local maximum. Since $(a_1,\ldots,a_5)\in E(H)$, we must have $\delta_1<\delta_2>\delta_3<\delta_4$ from (II) $m(e)=2$. Similarly, since $(a_2,\ldots,a_6)\in E(H)$ and $(a_3,\ldots,a_7)\in E(H)$, we have $\delta_4>\delta_5$ and $\delta_5<\delta_6$. Let $f=(a_1,\ldots,a_6)$, then we have $m(f)=3$. Note that $\delta_2$ is a local maximum, so $n(f\setminus a_4)=1$ from Claim \ref{mono}. Thus $m(f\setminus a_4)=1$, which implies that $f\setminus a_4\notin E(H)$ from the definition of the coloring, a contradiction.

Now we assume that $\delta_2$ is a local minimum. Since $(a_1,\ldots,a_5)\in E(H)$, we must have that $\delta_1>\delta_2<\delta_3>\delta_4$  from (II). Similarly, we have $\delta_4<\delta_5$ and $\delta_5>\delta_6$. Let $f=(a_2,\ldots,a_7)$, then we have $m(f)=3$. Note that $\delta_3$ is a local maximum, so $n(f\setminus a_4)=1$ from Claim \ref{mono}, which implies that $f\setminus a_4\notin E(H)$, again a contradiction.

Thus we have proved that $H$ is $K^{(5)}_7$-free. It remains to show that $\alpha_6(H)<4(\alpha_5(H'))^2$.

Set $n = 4t^2$ where $t=\alpha_5(H')$. On the contrary, there are vertices $a_1 < \cdots< a_n$ that
induce a $6$-independent set in $H$. If the corresponding sequence $\delta_1,\ldots,\delta_{n-1}$ contains fewer
than $4t$ local extrema, then there is some $j$ such that $\delta_j ,\ldots, \delta_{j+t}$ form a monotone sequence.
We may assume that this sequence is monotone decreasing.
Since $t = \alpha_5(H')$, these $t + 1$ vertices $\delta_j ,\ldots, \delta_{j+t}$ contain a copy of $K^{(4)}_5$ in $H'$. Say this copy is given by $\delta_{j_1},\ldots, \delta_{j_5}$. Then by Property F, the vertices $a_{j_1} < \cdots< a_{j_5} < a_{j_5+1}$ induce a copy
of $K^{(5)}_6$, a contradiction. Thus we may assume that the sequence $\delta_1,\ldots, \delta_{n-1}$ contains at least $4t$ local extrema.

Let $\Delta=\{\delta_{i_1},\ldots,\delta_{i_{4t}}\}$ be the first $4t$ local extrema, where $\delta_{i_j}=\delta(a_{i_j},a_{i_{j}+1})$. We shall find
a $6$-set $X$ in $H$ such that $X$ forms a $K^{(5)}_6$ in $H$. Note that at least $2t$ $\delta_{i_j}$'s among these local extrema are local maxima. Denote the first $2t$ local maxima by $\Delta_{\max}=\{\delta_{j_1},\ldots,\delta_{j_{2t}}\}$, then any two consecutive elements in $\Delta_{\max}$ are distinct from Property G. By a similar argument as above,  $\Delta_{\max}$ contains at least two local  extrema with one is local maximum, denoted by $\delta_{j_\ell}$, that is $\delta_{j_{\ell-1}}<\delta_{j_\ell}>\delta_{j_{\ell+1}}$. 

Let $\delta_{i_{p}}\in \Delta$ be the unique local minimum between $\delta_{j_{\ell-1}}$ and $\delta_{j_\ell}$; $\delta_{i_{q}}\in \Delta$ be the  unique local minimum between $\delta_{j_{\ell}}$ and $\delta_{j_{\ell+1}}$, respectively. Hence, $$\delta_{j_{\ell-1}}>\delta_{i_{p}}<\delta_{j_\ell}>\delta_{i_{q}}<\delta_{j_{\ell+1}}~~
\text{and}~~\delta_{j_{\ell-1}}<\delta_{j_\ell}>\delta_{j_{\ell+1}}.$$ 
Let $X=\{a_{j_{\ell-1}},a_{i_p},a_{i_{p}+1},a_{i_q},a_{i_{q}+1},a_{j_{\ell+1}+1}\}$.  By Property B, we have that $\delta(a_{j_{\ell-1}},a_{i_p})=\delta_{j_{\ell-1}}$, $\delta(a_{i_{p}+1},a_{i_q})=\delta_{j_\ell}$ and $\delta(a_{i_{q}+1},a_{j_{\ell+1}+1})=\delta_{j_{\ell+1}}$.
For convenience, rewrite $X=(a_1',\dots,a_6')$, and $\delta_i'=\delta(a_i',a'_{i+1})$.
Consider the $5$-tuple $e_i = X\setminus a'_i$ for any $1\le i\le 6$. We have that
\begin{align*}
\delta(e_1)&=(\delta_2'<\delta_3'>\delta_4'<\delta_5'),\\
\delta(e_6)&=(\delta_1'>\delta_2'<\delta_3'>\delta_4'),\\
\delta(e_2)=\delta(e_3)&=(\delta_1'<\delta_3'>\delta_4'<\delta_5'),\\
\delta(e_4)=\delta(e_5)&=(\delta_1'>\delta_2'<\delta_3'>\delta_5').
\end{align*}
 Thus, $X$ forms a $K^{(5)}_6$ in $H$ from (II), a contradiction.
This completes the proof.
\hfill$\Box$


\subsection{Proof of Theorem \ref{main-2}}

Since many details of the proof for $k=6$ are similar to that of $k\geq 7$, we shall prove the case for $k\geq 7$ at first.

\medskip\noindent
{\bf Part (I) \; $k\geq 7$.}

\medskip

We will show that $g(k,2^N)< (k+2) g(k-1, N)$ for $k\geq 7$. From the definition of $g(k-1, N)$, there is a $K_{k+1}^{(k-1)}$-free $(k-1)$-graph $H'$ on $\{0,1,\dots,N-1\}$ with $\alpha_{k}(H')=g(k-1, N)$. 
We shall construct a $K_{k+2}^{(k)}$-free $k$-graph $H$ on the vertex set $V(H)=\{0,1,\dots,2^N-1\}$
such that $\alpha_{k+1}(H)<(k+2) \alpha_{k}(H')$.



We define the edges of $H$ as follows. For any $k$-tuple $e=(a_1,a_2,\ldots,a_{k})$ of $V(H)$, set $e \in E(H)$ if and only if one of the following holds:

\medskip
(I) $m(e) =0$ and $\{\delta_1,\ldots,\delta_{k-1}\} \in E(H')$;

\smallskip
(II) $m(e) =k-3$;

\smallskip
(III)  $m(e) =k-4$ and $\de1<\de2$.

\medskip
From the above definition, we know $e\not\in E(H)$ if one of the following appears:

\smallskip
(IV) $m(e) =k-4$ and $\de1>\de2$;

\smallskip
(V)  $1\le m(e) \le k-5$.

\ignore{ In particular, if , or $1\le m(e) \le k-5$, then $e\not\in E(H)$. 
\medskip
(i) if $\delta_1,\ldots,\delta_{k-1 }$ form a monotone sequence, then $e \in E(H)$ $\Leftrightarrow$ $(\delta_1,\ldots,\delta_{k-1}) \in E(H')$;

\smallskip
(ii) if $\delta_1,\ldots,\delta_{k-1}$ is not monotone, then $e \in E(H)$  $\Leftrightarrow$  
$m(e) =k-3$, or $m(e) =k-4$ and $\de1<\de2$. 
}


\medskip
Moreover, Properties A, B and F hold (see in Section \ref{B-p}).


\medskip
We first show that $H$ contains no $(k+2)$-clique.
On the contrary, there exists a $(k+2)$-set $A=(a_1, \dots, a_{k + 2})$ that induces a $K_{k + 2}^{(k)}$ in $H$.  

\begin{claim} \label{four} Let $A=(a_1, \dots, a_{k + 2})$, which forms a $(k+2)$-clique in $H$ and $\dd(A)$ is nonmonotone, then there is no consecutive monotone subsequence of length $4$  of $\dd(A)$. Consequently, any consecutive subsequence of length $4$  of $\dd(A)$ contains a local extremum.
\end{claim}

\noindent
{\bf Proof.} On the contrary, $\de l,\de {l+1},\de {l+2},\de {l+3}$ form a consecutive monotone sequence of length $4$ of $\dd(A)$.
Since $\dd(A)$ is nonmonotone, we can extend this sequence to a new sequence $\de p,\ldots,\de {p+4}$ of length $5$ for some $p$ that contains a local extremum and two local monotone. Now we choose a $k$-tuple $e$ such that $e\supseteq [A]_{p}^{p+5}$. Thus 
$$1 \leq m(e)=(k-3)-n(e) \le (k-3)-2= k -5,$$ implying $e\not\in H$ from (V), a contradiction.
\hfill$\Box$

\medskip
 In the following, we separate the proof into four cases according to the number of local monotone and local extrema in $\delta(A)$.  ({\em Let us point out that all the cases but Case 3 will also lead to contradictions for $k=6$.})

\medskip
\emph{Case 1.}  $m(A)=0$, i.e., $\delta_2,\dots,\delta_{k}$ are local monotone.
For this case, $\delta_1,\ldots, \delta_{k + 1}$ form a monotone sequence.  By Property F, every $(k-1)$-tuple in the set $\{\delta_1,\ldots, \delta_{k + 1}\}$ forms an edge in $H'$ implying $H'$ contains a copy of $K_{k+1}^{(k-1)}$, which leads to a contradiction.


\medskip
\emph{Case 2.}  $m(A)\ge1$ and $n(A)\ge2$.    
Suppose first that $\delta_1 > \delta_2$. 
Let $e= A\setminus (a_{k+1},a_{k+2})$. Since $m(A)\ge1$, i.e., $\dd(A)$ is nonmonotone, we know from Claim \ref{four} that $\de {k-2},\de {k-1},\de {k},\de {k+1}$ cannot form a consecutive monotone sequence. So dropping $(a_{k+1},a_{k+2})$ from $A$ will destroy at most a local monotone, and thus $n(e)\ge1$. Moreover, $m(e)\ge1$ by noting $e\supseteq [A]_1^5$ and $\delta([A]_1^5)$ contains a local extremum from Claim \ref{four}. Hence $1\le m(e)\leq k-4.$ If $m(e)=k-4$, together with $\delta_1 > \delta_2$, then $e\notin H$ from (IV), a contradiction. If $1\le m(e)\le k-5$, then $e\notin H$ from (V), again a contradiction.

Now we assume $\delta_1 < \delta_2$. 
Suppose that $\de 2$ is a local monotone. Since $m(A)\geq1$, i.e., $\dd(A)$ is nonmonotone, then from Claim \ref{four}, we have $\delta_1 < \delta_2<\de3>\de4$; otherwise $\delta_1 < \delta_2<\de3<\de4$ would be a consecutive monotone subsequence of length $4$ in the sequence $\dd(A)$. Consider the $k$-tuple $e= A\setminus (a_1,a_{2})$. By deleting $(a_1,a_{2})$, it destroys the local monotone $\de 2$. So   $n(e)\ge1$ by noting $n(A)\ge2$. 
Moreover, $m(e)\ge1$ since $\delta([A]_3^7)$ contains a local extremum from Claim \ref{four}. Hence $1\le m(e)\leq k-4.$ Note that $\delta_3 > \delta_4$. By a similar argument as above, we have $e\notin H$ from (IV) or (V), a contradiction.
Therefore, we may assume that $\de 2$ is a local extremum, and so $\delta_1 <\de2>\de3$. Now consider the $k$-tuple $e= A\setminus (a_1,a_{k+2})$.
Since dropping $a_{k+2}$ from $A$ will destroy at most a local monotone, we have $n(e)\geq1$. Moreover, $m(e)\ge1$ by noting $\delta([A]_2^6)$ contains a local extremum from Claim \ref{four}. Hence $1\le m(e)\leq k-4.$ Together with $\delta_2 > \delta_3$ we know $e\notin H$ by (IV) or (V),  again a contradiction.

\medskip
\emph{Case 3.} $n(A)=1$.    
Suppose first $\delta_1 > \delta_2$. 
Let $\de l$ be the unique local monotone in $\dd(A)$. Since $(k+2)-3-4\geq2$, we can choose a $k$-tuple $e= A\setminus (a_p,a_q)$ where $a_p\neq a_q$ such that $a_p, a_q\notin [A]_1^3\cup[A]_{l-1}^{l+2}$. Since all $\delta_i$ for $2\le i\le k$ but $i=l$ are local extrema, the number of local monotone of the resulting sequence of $\delta$'s by deleting $(a_p,a_q)$ may increase by at most 2. Therefore, $1\le n(e)\le3$ implying $1\le m(e)\le k-4$. Together with $\delta_1 > \delta_2$ we know $e\notin H$ by (IV) or (V),  a contradiction.

Now we assume  $\delta_1 < \delta_2$. 
If $\de 2$ is the unique local monotone, then $m([A]_2^7)=|\{\delta_3,\delta_4,\delta_5\}|=3$ and $\de 3$ is a local maximum.
So we apply Claim \ref{mono} to $[A]_2^7$ to obtain $n([A]_2^7\setminus a_4)=1$ and $\delta_3$ or $\delta_5$ would be a local monotone of $[A]_2^7\setminus a_4$. Consider the $k$-tuple $e= A\setminus (a_4,a_{k+2})$. Thus $n(e)=2$, implying $m(e)=k-5$. By (V), we know $e\notin H$, a contradiction.
Therefore, we may assume that $\de 2$ is a local extremum.  Then $\de1<\de2>\de3$.
Since $(k+2)-4-4\geq1$, we can choose a $k$-tuple $e= A\setminus (a_1,a_p)$ such that $a_p\notin [A]_1^4\cup[A]_{l-1}^{l+2}$. Since all $\delta_i$ where $2\le i\le k$ but $i=l$ are local extrema, the number of local monotone of the resulting sequence of $\delta$'s by deleting $(a_1,a_p)$ may increase by at most 1. Therefore, $1\le n(e)\le2$ implying $k-5\le m(e)\le k-4$. Together with $\delta_2 > \delta_3$ we know $e\notin H$ by (IV) or (V),  a contradiction.

\medskip
\emph{Case 4.}  $n(A)=0$.   
Suppose first  $\delta_1 > \delta_2$.
Since $m([A]_2^7)=|\{\delta_3,\delta_4,\delta_5\}|=3$ and $\de 3$ is a local maximum,
we apply Claim \ref{mono} to $[A]_2^7$ to obtain that $n([A]_2^7\setminus a_4)=1$. Consider the $k$-tuple $e= A\setminus (a_4, a_{k+2})$.
Then $n(e)=1$ implying $m(e)= k-4$. Together with $\delta_1 > \delta_2$, we know $e\notin H$ by (IV), a contradiction.
Now we assume 
$\delta_1 < \delta_2$. Then $\delta_2 > \delta_3$ by noting $n(A)=0$.
Since $m([A]_3^8)=|\{\delta_4,\delta_5,\delta_6\}|=3$ and $\de 4$ is a local maximum, we apply Claim \ref{mono} to $[A]_3^8$ to obtain that $n([A]_3^8\setminus a_5)=1$.
Consider the $k$-tuple $e= A\setminus (a_1, a_5)$.
A similar argument as above yields that $e\notin H$ by (IV), again a contradiction.

Thus we have proved that $H$ is $K_{k + 2}^{(k)}$-free. It remains to show $\alpha_{k+1}(H)<(k+2)\alpha_k(H').$ 

Set $n = (k+2)t$ where $t = \alpha_k(H')$. On the contrary, there are vertices $a_1< \cdots< a_n$ that induce a $(k+1)$-independent set in $H$. If the corresponding sequence $\delta_1,\ldots, \delta_{n-1}$ contains fewer than $k+2$ local extrema, then there is a $j$ such that $\delta_j, \ldots, \delta_{j+t}$ form a monotone sequence. Without loss of generality, we may assume that it is a monotone decreasing sequence. Since $t = \alpha_k(H')$, these $t + 1$ vertices $\delta_j,\ldots,\delta_{j + t}$ contain a copy of $K^{(k-1)}_k$ in $H'$.  Say this copy is given by $\delta_{j_1},\ldots, \delta_{j_k}$.  Then by Property F, the vertices $a_{j_1} < \cdots < a_{j_k} < a_{j_k + 1}$ induce a copy of $K_{k + 1}^{(k)}$, a contradiction. Thus we may assume that the sequence $\delta_1,\ldots ,\delta_{n-1}$ contains at least $k+2$ local extrema. 

Let $\delta_{i_1},\ldots, \delta_{i_{k+2}}$ be the first $k+2$ local extrema, where $\delta_{i_j} = \delta(a_{i_j},a_{i_j + 1})$. We shall find a $(k+1)$-set $X$ in $H$ such that the number of local extrema $m(X)=k-2$.   
Note that at least $\lceil(k+1)/2\rceil$ $\delta_{i_j}$'s among these local extrema are local minima, 
 and we let $X$ include the first $k+1$ vertices out of all $(a_{i_j},a_{i_j + 1})$ where $\delta_{i_j}$ is a local minimum. For convenience, we rewrite $X=(a'_1,\ldots, a'_{k+1})$ and $\delta'_i = \delta(a'_i,a'_{i + 1})$.
By Property B,  we have that $\delta'_{2i}$ is a local maximum for $1\le i\le\lfloor(k-1)/2\rfloor$, and $\delta'_{2i+1}$ is a local minimum for $1\le i\le\lfloor(k-2)/2\rfloor$. Namely,
 \[\delta'_1<\delta'_2>\delta'_3<\cdots \delta'_k.\]
Therefore, the sequence $\{\delta'_1,\ldots, \delta'_k\}$ has $k-2$ local extrema, i.e., $m(X)=k-2$ as desired. 

Now we consider the $k$-tuple $e_i = X\setminus a'_i$ for $1\le i\le k+1$. We claim $e_i\in H$ for $1\le i\le k+1$.
Since $m(e_1)=k-3$, we have $e_1\in H$  from (II). Similarly, we have  $e_2\in H$ by noting that $\delta(a_1',a_3')=\delta_2'$ from Property B.
For each $3\leq i\leq k+1$, $e_i$ satisfies $\de 1'<\de2'$. 
Note that dropping $a_i'$ for $3\leq i\leq k+1$ from $X$ will produce at most one local monotone. 
Thus for $3\leq i\leq k+1$, $0\leq n(e_i)\leq 1$ and so $k-4 \leq m(e_i)\leq k-3$.
Together with $\de 1'<\de2'$, we have $e_i\in H$ for $3\leq i\leq k+1$ from (II) or (III).
Therefore, $X$ induces a $K_{k + 1}^{(k)}$ in $H$, a contradiction.
\hfill$\Box$

\bigskip\noindent
{\bf Part (II) \; $k=6$.}

\medskip

We will show $g(6,2^N)< 8 g(5, N)$. From the definition of $g(5, N)$, there is a $K_{7}^{(5)}$-free $5$-graph $H'$ on $\{0,1,\dots,N-1\}$ with $\alpha_{6}(H')=g(5, N)$. 
We shall construct a $K_{8}^{(6)}$-free $6$-graph $H$ on the vertex set $V(H)=\{0,1,\dots,2^N-1\}$ 
such that $\alpha_{7}(H)<8 \cdot\alpha_{6}(H')$.

We define the edges of $H$ as follows. For any $6$-tuple $e=(a_1,a_2,\ldots,a_{6})$ of $V(H)$, set $e \in E(H)$ if and only if one of the following holds:

\medskip
(I) $m(e) =0$ and $\{\delta_1,\ldots,\delta_{5}\} \in E(H')$;

\smallskip
(II)  $m(e)\in\{2,3\}$ and $\de1<\de2$;

\smallskip
(III) $m(e)=3$, $\de1>\de2$ and  either $\de1>\de3$ or $\de3<\de5$.

\medskip
From the above definition, we know $e\not\in E(H)$ if one of the following appears:

\smallskip
(IV)  $m(e) =1$; 

\smallskip
(V) $m(e) =2$ and $\de1>\de2$;

\smallskip
(VI) $m(e)=3$, $\de1>\de2$ and  $\de1<\de3>\de5$.

\ignore{ In particular, if , or $1\le m(e) \le k-5$, then $e\not\in E(H)$. 
\medskip
(i) if $\delta_1,\ldots,\delta_{k-1 }$ form a monotone sequence, then $e \in E(H)$ $\Leftrightarrow$ $(\delta_1,\ldots,\delta_{k-1}) \in E(H')$;

\smallskip
(ii) if $\delta_1,\ldots,\delta_{k-1}$ is not monotone, then $e \in E(H)$  $\Leftrightarrow$  
$m(e) =k-3$, or $m(e) =k-4$ and $\de1<\de2$. 
}


\medskip
Moreover, Properties A, B, D, F and G hold (see in Section \ref{B-p}).

\begin{claim} \label{mono-2} 
Given $f = (a_1,\ldots,  a_5)$, if $\de 1<\de 2>\de 3<\de4$ and $\de2<\de4$, then $\de 2$ is a local monotone of $f\setminus a_3$; if $\de 1>\de 2<\de 3>\de4$ and $\de1>\de3$, then $\de 3$ is a local monotone of $f\setminus a_3$. 
\end{claim}

\noindent
{\bf Proof.}  By Property B, for the first case, we have $\dd(f\setminus a_3) =(\de1<\de2<\de4)$; for the second case, $\dd(f\setminus a_3) =(\de1 >\de3>\de4)$. The assertion follows.
\hfill$\Box$

\medskip
We first show that $H$ contains no $8$-clique.
On the contrary, there is a set $A=(a_1, \dots, a_{8})$ that induces a $K_{8}^{(6)}$ in $H$.  
We will find a $6$-tuple $e\subseteq A$ with $m(e)=1$; or $m(e)=2$ and $\de1>\de2$; or $m(e)=3$,  $\de1>\de2$ and $\de1<\de3>\de 5$.
Hence $e\notin H$ from (IV), (V), or (VI), a contradiction. 
The discussion is the same as in Part (I) for all cases except the case that $n(A)=1$. Therefore, we may assume that $n(A)=1$ in the following.  

\medskip
\emph{Case 1.}  $\delta_1 > \delta_2$. 
If $\de {i}$ is the unique local monotone where $2\leq i\leq 4$, then we consider the $6$-tuple $e= A\setminus (a_7,a_8)$. 
Clearly $n(e)=1$, implying $m(e)=2$. Together with $\delta_1 > \delta_2$, we have $e\notin H$ from (V), a contradiction. 
If $\de 6$ is the unique local monotone, then we consider the $6$-tuple $e= A\setminus (a_1,a_4)$. Note that $m([A]_2^7)=|\{\delta_3,\delta_4,\delta_5\}|=3$ and $\de 3$ is a local maximum. So we apply Claim \ref{mono} to $[A]_2^7$ to obtain that $n([A]_2^7\setminus a_4)=1$. Thus $n(e)=2$, implying $m(e)=1$. Hence $e\notin H$ from (IV), a contradiction.

It remains to verify the case where $\de 5$ is the unique local monotone. 
Then $\delta_1 > \delta_2<\de3>\de4$, and so we have $\de3\neq\de 5$ by Property D.
If $\de 3<\de 5$, then we consider the $6$-tuple $e= A\setminus (a_4,a_8)$. 
Since $\de 2<\de 3>\de 4<\de5$ and $\de3<\de5$, we apply Claim \ref{mono-2} to $[A]_2^6$ to obtain that $\de3$ is a local monotone of $[A]_2^6\setminus a_4$. Together with $\de 5$ is a local monotone of $e$, we have $n(e)=2$, implying $m(e)=1$. Thus $e\notin H$ from (IV), a contradiction.
If $\de 3>\de 5$, then we consider the $6$-tuple $e= A\setminus (a_4,a_5)$. Note that $\dd (e)=(\de1>\de2<\de3,\;\de6>\de7)$, and $\delta(a_3,a_5)=\delta_3$ by Property B. 
Since $\de3>\de5$, we apply Property D to $(a_3,a_5,a_6,a_7)$ to obtain that $\de3\neq\de 6$. 
If $\de3>\de6$, then $\de 6$ is a local monotone; otherwise $\de 3$ is a local monotone.
Thus $n(e)=1$, implying $m(e)=2$. Together with $\delta_1 > \delta_2$ we know $e\notin H$ by (V),  a contradiction.

\medskip
\emph{Case 2.}  $\delta_1 < \delta_2$. 
If $\de 2$ is the unique local monotone, then $m([A]_2^7)=|\{\delta_3,\delta_4,\delta_5\}|=3$ and $\de 3$ is a local maximum.
So we apply Claim \ref{mono} to $[A]_2^7$ to obtain that $n([A]_2^7\setminus a_5)=1$ and $\delta_3$ or $\delta_5$ would be a local monotone of $[A]_2^7\setminus a_5$.
Consider the $6$-tuple $e= A\setminus (a_5,a_{8})$. Together with that $\delta_2$ is a local monotone, we have $n(e)=2$ implying $m(e)=1$. By (IV), we know $e\notin H$, a contradiction.
If $\de {i}$ is the unique local monotone where $3\leq i\leq 5$. For this case, we have $\delta_1 < \delta_2 > \delta_3$. Consider the $6$-tuple $e= A\setminus (a_1,a_{8})$. We have $n(e)=1$, implying $m(e)=2$. Together with $\delta_2 > \delta_3$ we know $e\notin H$ by (V), a contradiction.

It remains to verify the case where  $\de {6}$ is the unique local monotone. Consider the $6$-tuple $e= A\setminus (a_1,a_{8})$. 
Note that $\de{2i}>\de{2i+1}$ for $i=1,2$, by Property D we have that $\de2\neq\de 4$ and $\de 4\neq\de 6$.
We claim that $\de 2<\de 4>\de 6$.
Otherwise, if $\de 2>\de 4$, then we consider the $6$-tuple $e= A\setminus (a_1,a_4)$. 
Since $\de 2>\de 3<\de 4>\de5$ and $\de2>\de4$, we apply Claim \ref{mono-2} to $[A]_2^6$ to obtain that $\delta_4$ is a local monotone of $[A]_2^6\setminus a_4$. Thus $\de 4$ and  $\de 6$ are two local monotone of $e$, i.e., $n(e)=2$, implying $m(e)=1$. Hence $e\notin H$ from (IV), a contradiction. 
If $\de 4<\de 6$, then we consider the $6$-tuple $e= A\setminus (a_1,a_{5})$.
Since $\de 3<\de 4>\de 5<\de6$ and $\de4<\de6$, we apply Claim \ref{mono-2} to $[A]_3^7$ to obtain that $\delta_4$ is a local monotone of $[A]_3^7\setminus a_5$. Since $\de 4<\de 6<\de 7$, we have that $\de6$ is still a local monotone of $e$. Thus $n(e)=2$, implying $m(e)=1$. By (IV), we know $e\notin H$, a contradiction.  
Thus $6$-tuple $e= A\setminus (a_1,a_{8})$ satisfies $\de 2>\de3$ and $\de 2<\de 4>\de 6$, and $m(e)=|\{\delta_3,\delta_4,\delta_5\}|=3$, so by (VI), $e\notin H$, a contradiction. 

Thus we have proved that $H$ is $K_{8}^{(6)}$-free. It remains to show $\alpha_{7}(H)<8\cdot\alpha_6(H').$ 

Set $n = 8t$ where $t = \alpha_6(H')$. On the contrary, there are vertices $a_1< \cdots< a_n$ that induce a $7$-independent set in $H$. If the corresponding sequence $\delta_1,\ldots, \delta_{n-1}$ contains fewer than $8$ local extrema, then by a similar argument as Part (I), we can find a copy of $K_{7}^{(6)}$ in $H$, a contradiction. Thus we may assume that the sequence $\delta_1,\ldots ,\delta_{n-1}$ contains at least $8$ local extrema. 
Therefore the sequence $\delta_1,\ldots ,\delta_{n-1}$ contains at least $4$ local maxima. 

We write $M$ for the sequence consisting of the first $4$ local maxima.
By Property G, we know that any two consecutive local maxima of these $4$ local maxima are distinct.  
Therefore, we can find three consecutive local maxima $\de {k_1},\de {k_2},\de {k_3}$ which form a monotone sequence or $\de {k_1}>\de {k_2}<\de {k_3}$, where $\de {k_i}=\dd (a_{k_i},a_{{k_i}+1})$ for $1\leq i\leq 3$.

For $i=2,3$, let $\de {l_i}=\delta(a_{l_i},a_{l_i+1})$ be the closest local minimum on the left of $\de {k_i}$, and let 
\[
X=(a_{k_1-1},a_{k_1},a_{l_2},a_{l_2+1},a_{l_3},a_{l_3+1},a_{k_3+1}).
\]
By Property B, $\delta(a_{k_1},a_{l_{2}})=\delta_{k_1}$, $\delta(a_{l_2+1},a_{l_3})=\delta_{k_2}$, and $\delta(a_{l_3+1},a_{k_3+1})=\delta_{k_3}$.
We obtain that the corresponding sequence satisfies that $$\de {k_1-1}<\de {k_1}>\de {l_2}<\de {k_2}>\de {l_3}<\de {k_3},$$ which contains $4$ local extrema, i.e., it has no local monotone. 

In the following, we shall show that each $6$-tuple of $X$ is an edge in $H$.
For simplicity, rewrite $X=(a_1',\dots,a_7')$, and $\delta_i'=\delta(a_i',a'_{i+1})$.
Consider the $6$-tuple $e_i = X\setminus a'_i$ for $1\le i\le 7$. 
By Property B,  we have that
\begin{align*}
\delta(e_1)=\delta(e_2)&=(\delta_2'>\delta_3'<\delta_4'>\delta_5'<\delta_6'),\\
\delta(e_3)=\delta(e_4)&=(\delta_1'<\delta_2', \;\delta_4'>\delta_5'<\delta_6'),\\
\delta(e_5)=\delta(e_6)&=(\delta_1'<\delta_2'>\delta_3'<\delta_4', \; \delta_6'),\\
\delta(e_7)&=(\delta_1'<\delta_2'>\delta_3'<\delta_4'>\delta_5').
\end{align*}
For  $e_1$ and $e_2$, we know that $m(e_1)=m(e_2)=3$ and  $\de 2'>\de3'$.
Moreover, since $\delta'_2,\delta'_4,\delta'_6$ are monotone or $\delta'_2>\delta'_4<\delta'_6$ from the above, then either $\delta'_2>\delta'_4$ or $\delta'_4<\delta'_6$. Therefore, $e_1,e_2\in H$ from (III). For $e_3$ and $e_4$, either $\delta_2'>\delta_4'$ or $\delta_2'<\delta_4'$, thus we have $m(e_3)=m(e_4)=2$. Together with $\delta_1'<\delta_2'$ yields that $e_3, e_4\in H$ from (II). For $e_i$ where $5\le i\le 7$, we have $m(e_i)=2$ or $3$. Together with $\delta_1'<\delta_2'$ yields that $e_i\in H$ for $5\le i\le 7$ from (II) again.

Therefore, $X$ induces a $K_{7}^{(6)}$ in $H$. This final contradiction completes the proof.\hfill$\Box$

\end{spacing}
\end{document}